\documentclass[reqno,a4paper,12pt]{amsart}

\usepackage{enumitem}
\setenumerate{label=\textnormal{(\arabic*)}}

\usepackage{amsmath,amssymb,dsfont,verbatim,bm,array,mathtools}
\usepackage[latin1]{inputenc}
\usepackage{booktabs}
\usepackage[raggedright]{titlesec}
\usepackage{mathtools}
\usepackage{graphicx}
\usepackage{fancyvrb}
\usepackage[dvipsnames]{xcolor}

\newenvironment{blockcode}
  {\leavevmode\small\color{blue}\verbatim}
  {\endverbatim}

\newenvironment{figurecode}
  {\leavevmode\small\color{Sepia}\verbatim}
  {\endverbatim}

\titleformat{\chapter}[display]
{\normalfont\huge\bfseries}{\chaptertitlename\\thechapter}{20pt}{\Huge}
\titleformat{\section}
{\normalfont\Large\bfseries\center}{\thesection}{1em}{}
\titleformat{\subsection}
{\normalfont\large\bfseries}{\thesubsection}{1em}{}
\titleformat{\subsubsection}[runin]
{\normalfont\normalsize\bfseries}{\thesubsubsection}{1em}{}
\titleformat{\paragraph}[runin]
{\normalfont\normalsize\bfseries}{\theparagraph}{1em}{}
\titleformat{\subparagraph}[runin]
{\normalfont\normalsize\bfseries}{\thesubparagraph}{1em}{}
\titlespacing*{\chapter} {0pt}{50pt}{40pt}
\titlespacing*{\section} {0pt}{3.5ex plus 1ex minus .2ex}{2.3ex plus .2ex}
\titlespacing*{\subsection} {0pt}{3.25ex plus 1ex minus .2ex}{1.5ex plus .2ex}
\titlespacing*{\subsubsection}{0pt}{3.25ex plus 1ex minus .2ex}{1.5ex plus .2ex}
\titlespacing*{\paragraph} {0pt}{3.25ex plus 1ex minus .2ex}{1em}
\titlespacing*{\subparagraph} {\parindent}{3.25ex plus 1ex minus .2ex}{1em}

\input xypic
\xyoption{all}

\begin{document}
\title{Visually pleasing knot projections in R}

\author{Robin Hankin}
\address{AUT}
\email{hankin.robin@gmail.com}

\begin{abstract}
In this short article I introduce the {\tt knotR}
package~\cite{hankin2016}, which creates two dimensional knot diagrams
optimized for visual appearance using the R programming
language~\cite{rcore2016}.  The {\tt knotR} package is a systematic
R-centric suite of software for the creation of production-quality
artwork of knot diagrams, released under GPL2.
\end{abstract}

\maketitle

\section{Introduction}

A {\em mathematical knot} is a smooth, unoriented embedding of a
circle~$\mathbb{S}^1$ into~$\mathbb{R}^3$~\cite{manturov2004}.  Two
knots are said to be equivalent if one can be continuously deformed in
to the other; if so, there is a
homeomorphism~$h\colon\mathbb{R}^3\longrightarrow\mathbb{R}^3$ which
takes one embedded circle to the other.

\begin{figure}[h]
  \centering
  \includegraphics[width=10cm]{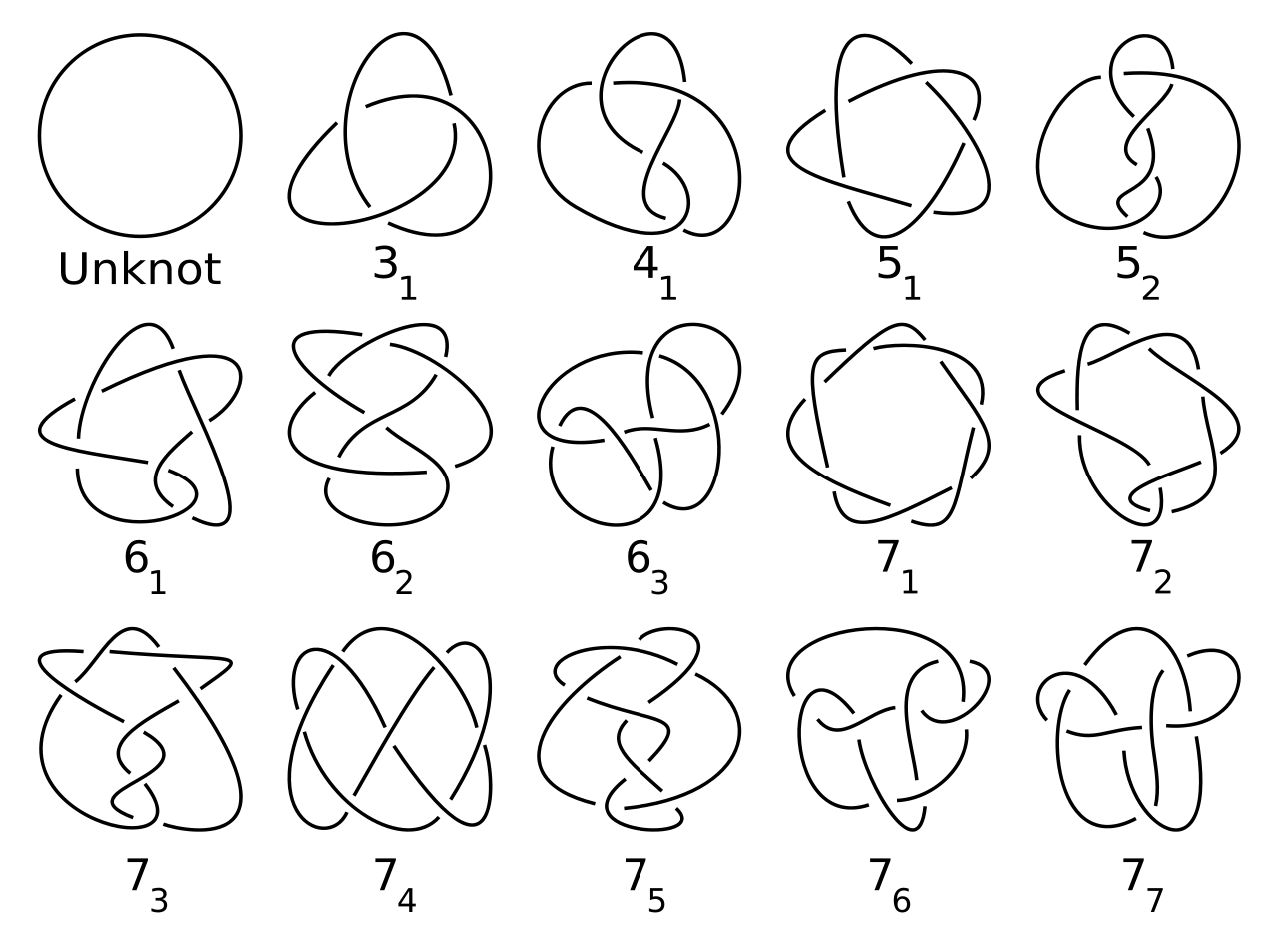}
  \caption{A table of prime knots up \label{knot_table} to seven
    crossings, labels following~\cite{alexander1926}.  Image taken
    from~\cite{wikipedia_knot_theory}.}
\end{figure}

It is common to present a knot using diagrams such as
Figure~\ref{knot_table}, in which a two-dimensional projection of the
knot is given with broken lines indicating where one strand passes
over another. 

Consider Figure~\ref{knot_table} from an \ae sthetic perspective; the
diagrams are representative of those available under a free license.
However, these diagrams are not suitable for high-quality artwork such
as posters: they are not vectorized.  Many of the underlying knots
possess a line of symmetry (at least, the diagrams do if the breaks
are ignored), which is not present in the visual representation.
Also, several of the strands cross at acute angles.  The diagram for
knot~$7_3$, for example, contains ugly kinks and abrupt changes of
curvature.

Such considerations suggest that knot diagrams might be produced by
minimization of some objective function that quantifies the visual
inelegance of a knot diagram.  The precise nature of such an objective
function is, of course, a subjective choice but one might require the
following desiderata:

\begin{itemize}
\item Curvature to be as uniform as possible
\item Strands to cross at right angles
\item Crossing points to be separate from one another
\item Any symmetry present in the knot should be enforced exactly, and be visually apparent
\end{itemize}

Knot diagrams may be created using vectorized graphics software such
as inkscape~\cite{kirsanov2009}: one specifies a sequence of control
points, then interpolates between these points to create a knot
diagram with the appropriate topology (inkscape is a widely-used
system available under the GPL).  One way of smoothly interpolating
between specified points is Bezier curves~\cite{olsen2014}.  A Bezier
curve is a visually pleasing polynomial path that can be used to
specify the path of a knot projection.

The package presented here allows one to specify a knot in terms of
its Bezier control points within inkscape, import the object into
R~\cite{rcore2016}, and then to use numerical optimization techniques
to improve the visual appearance of the knot.

\section{The package in use}

In this section, I give workflows for creating two simple knots:
firstly~$7_6$, followed by the figure-of-eight knot~$4_1$ which
requires imposition of symmetry constraints.

\begin{figure}[h]
  \centering
    \includegraphics[width=10cm]{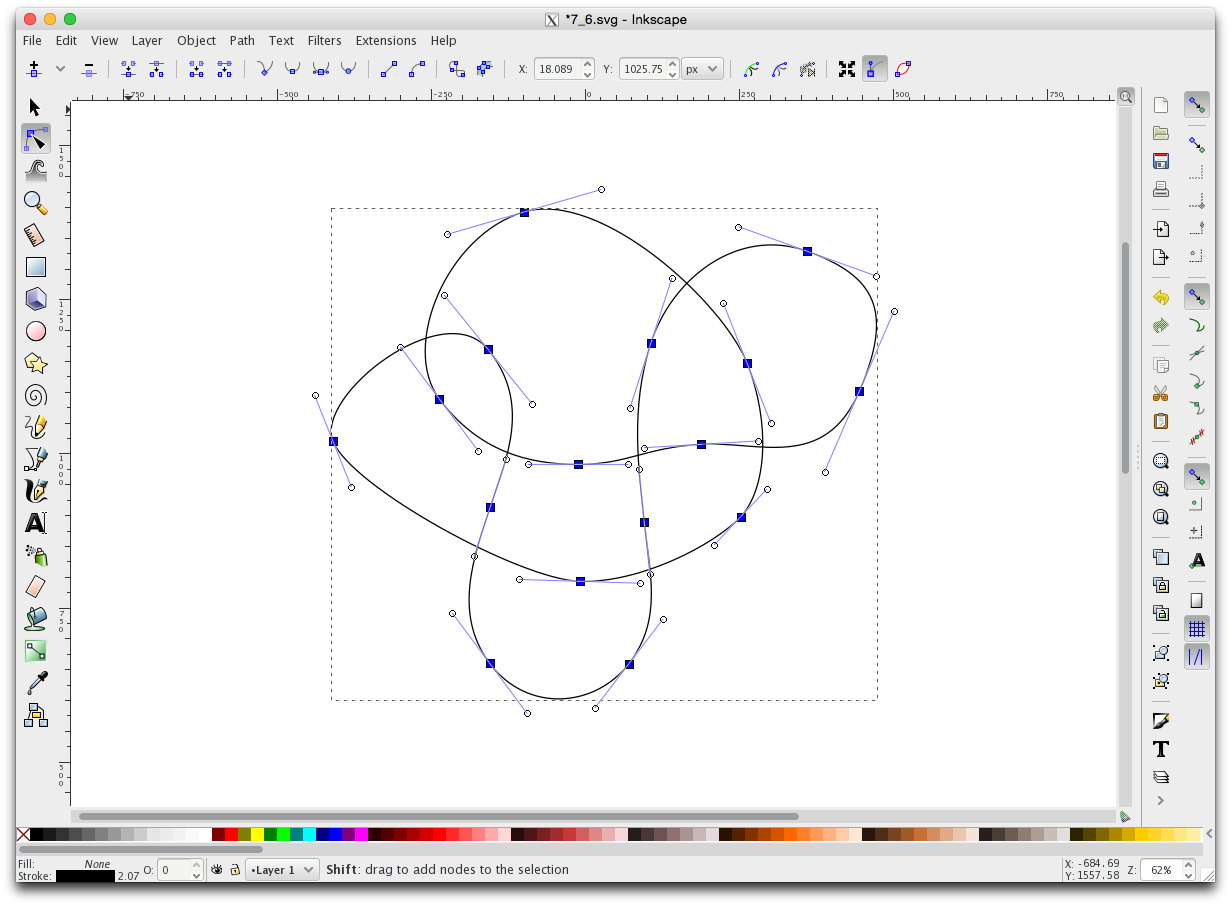}
\caption{Screenshot of inkscape\label{screenshot_inkscape_7_6} setup
  for knot~$7_6$.  Nodes are shown as squares, handles as small
  circles, symmetrically placed on either side of nodes}
\end{figure}

The first step is to create a closed curve in inkscape that shows the
rough outline of the knot~(Figure~\ref{screenshot_inkscape_7_6} shows
a screenshot of {\tt 7\_6\_first\_draft.svg}, supplied with the
package).  Note that this file contains only the knot {\em path}; the
over and under information is to be added later.

Here, knot paths are required to have Bezier handles that are
symmetrically placed with regard to nodes.  Although radius of
curvature is not necessarily matched at either side of a node, visual
continuity of the strand is ensured.  In the diagram, strand crossing
points are far from nodes, as this ensures visual continuity of
strands at crossing points, especially the understrand.

The knot shown in Figure~\ref{screenshot_inkscape_7_6} is clearly
suboptimal: even though the nodes are connected by Bezier curves which
are individually smooth, the path as a whole is ugly and unattractive
as its path has unsightly regions where the radius of curvature
changes abruptly.

Although it is possible in principle to improve the visual appearance
of the knot path by hand in inkscape, this is a surprisingly difficult
and frustrating task.  In order to remedy the flaws of the diagram
using an automated system, we first read the {\tt .svg} file into R
using the {\tt reader()} function; a typical R session follows:

\begin{blockcode}
> library(knotR)
> k76 <- reader(system.file("7_6_first_draft.svg",package="knotR"))
> head(k76)

              x         y
[1,]  -98.81963 339.81898
[2,] -223.87754 303.35366
[3,] -299.84521 121.06064
[4,] -236.36319  36.35340
[5,] -172.88118 -48.35384
[6,]  -92.86186 -69.78212
\end{blockcode}

Object~{\tt k76} is stored as an object of class {\tt inkscape}: a
two-column matrix with rows corresponding to the node and handle
positions of the inkscape path.  This representation has a certain
amount of redundancy as knot paths have handles which are
symmetrically placed with respect to nodes; also, the first node is
the same as the last for the loop is closed.  The package can coerce
inkscape objects to other forms, specifically {\tt minobj} objects,
which contain no redundancy (the position of each node, as well as one
of the handles, is stored); or {\tt controlpoints} objects, which
allow for easy construction of Bezier interpolation between nodes.

\begin{figure}[htbp]
  \begin{center}
\begin{figurecode}
> k76_rough <- reader(system.file("7_6_first_draft.svg",package="knotR"))
> knotplot2(k76_rough, seg=TRUE)
\end{figurecode}
\begin{centering}
\includegraphics{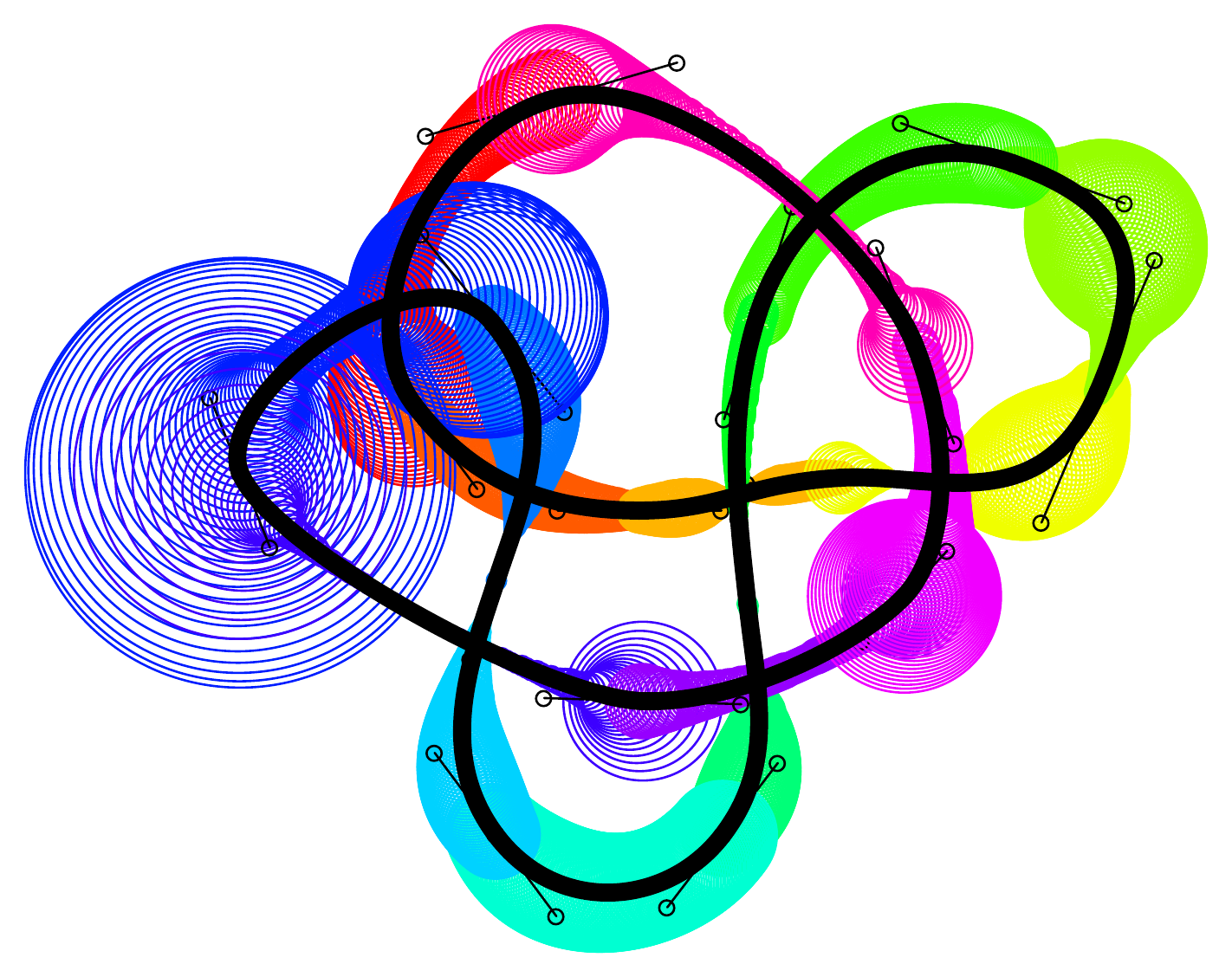}
\caption{The \emph{path} of (unoptimized) \label{7_6_rough}
  knot~$7_6$, showing Bezier handles as thin straight lines and
  circles.  The coloured circles have a radius  proportional to the
  curvature (that is, the reciprocal of the radius of curvature) along
  the strand; note large curvature at loop on left}
  \end{centering}
  \end{center}
\end{figure}

\begin{figure}[!tbp]
\begin{figurecode}
knotplot2(k7_6,text=TRUE,lwd=1,circ=0,rainbow=TRUE)
\end{figurecode}
 \centering
\includegraphics{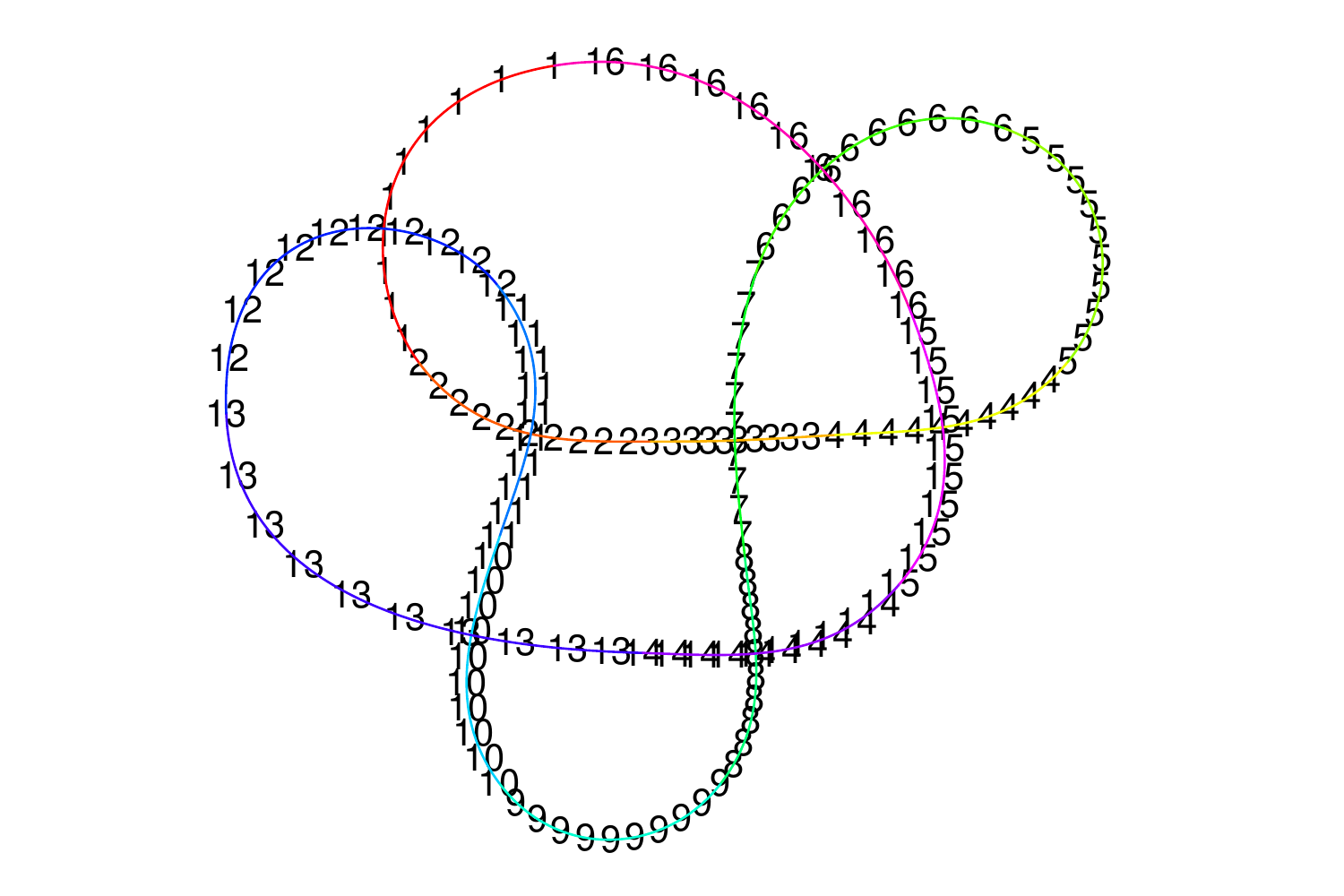}
\caption{Knot~$7_6$ with strands numbered \label{k76_strands} to guide
creation of over and under information} \end{figure}

To beautify it we need to specify a function of the path that
quantifies its displeasingness, and then minimize this objective
function using numerical methods.

Two examples of desiderata for such an objective function might be to
keep the strands crossing at right angles, and the overall bending
energy.  These are evaluated in the package by functions
{\tt total\_crossing\_angles()} and~{\tt total\_bending\_energy()}
respectively:

\begin{blockcode}
> b <- as.controlpoints(k76_rough)
> total_crossing_angles(b)

[1] 0.3145033

> total_bending_energy(b)

[1] 0.1276257
\end{blockcode}

The knots supplied in the package minimize a weighted sum of these and
other badnesses\footnote{Function {\tt badness()} includes various
``housekeeping'' badnesses which are used to make sure that the
minimum found by {\tt nlm()} is topologically identical to the
starting configuration.  Function
{\tt non\_crossing\_strand\_close\_approach\_badness()}, for example,
makes non-crossing strands ``repel'' one another so as not to
introduce spurious intersections.}, as evaluated by function
{\tt badness()}:

\begin{blockcode}
> badness(k76_rough)

[1] 6.16279
\end{blockcode}

This function may be minimized by numerical optimization:

\begin{blockcode}
> o <- nlm(badness, as.knotvec(k76_rough))
> k7_6 <- as.minobj(o$estimate)
> badness(k7_6)

  [1] 3.550152
\end{blockcode}

\begin{figure}[!tbp]
\begin{figurecode}
knotplot2(k7_6)
\end{figurecode}
  \centering
\includegraphics{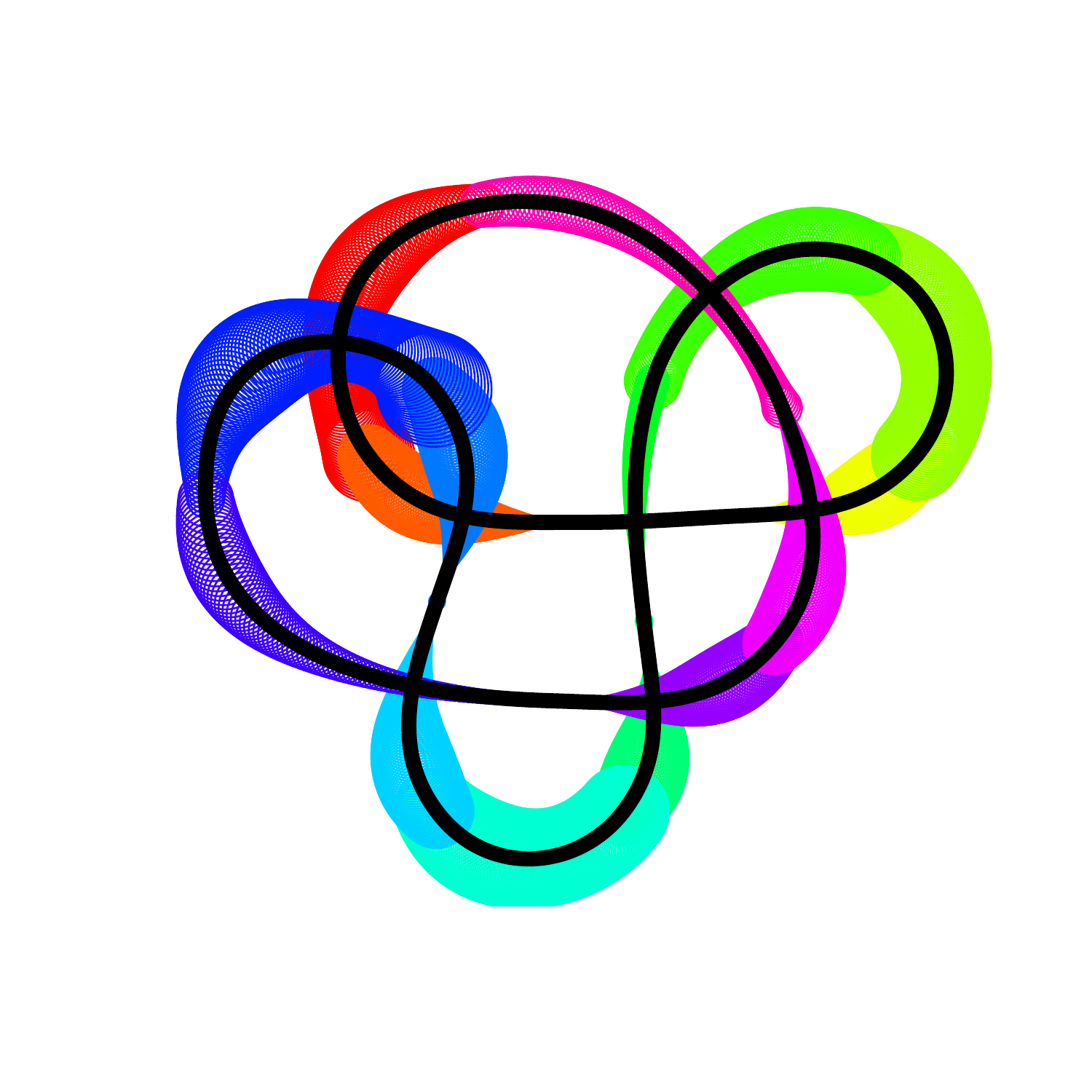}
\caption{Knot $7_6$, post-optimization\label{7_6}}
\end{figure}

(the above takes about an hour of CPU time: it is optimizing a
function of~64 real variables, and the objective function takes a few
seconds to evaluate).  However, the result is much nicer
(Figure~\ref{7_6}).

To specify the senses of the knot's\ 
crossings, we create an overunder object which is a two-column matrix:

\begin{blockcode}
> ou76 <- matrix(c(
+     12,01,
+     02,11,
+     07,03,
+     04,15,
+     16,06,
+     14,08,
+     10,13
+     ),byrow=TRUE,ncol=2)
\end{blockcode}

\begin{figure}[!tbp]
\begin{figurecode}
knotplot(k7_6)
\end{figurecode}
  \centering
\includegraphics{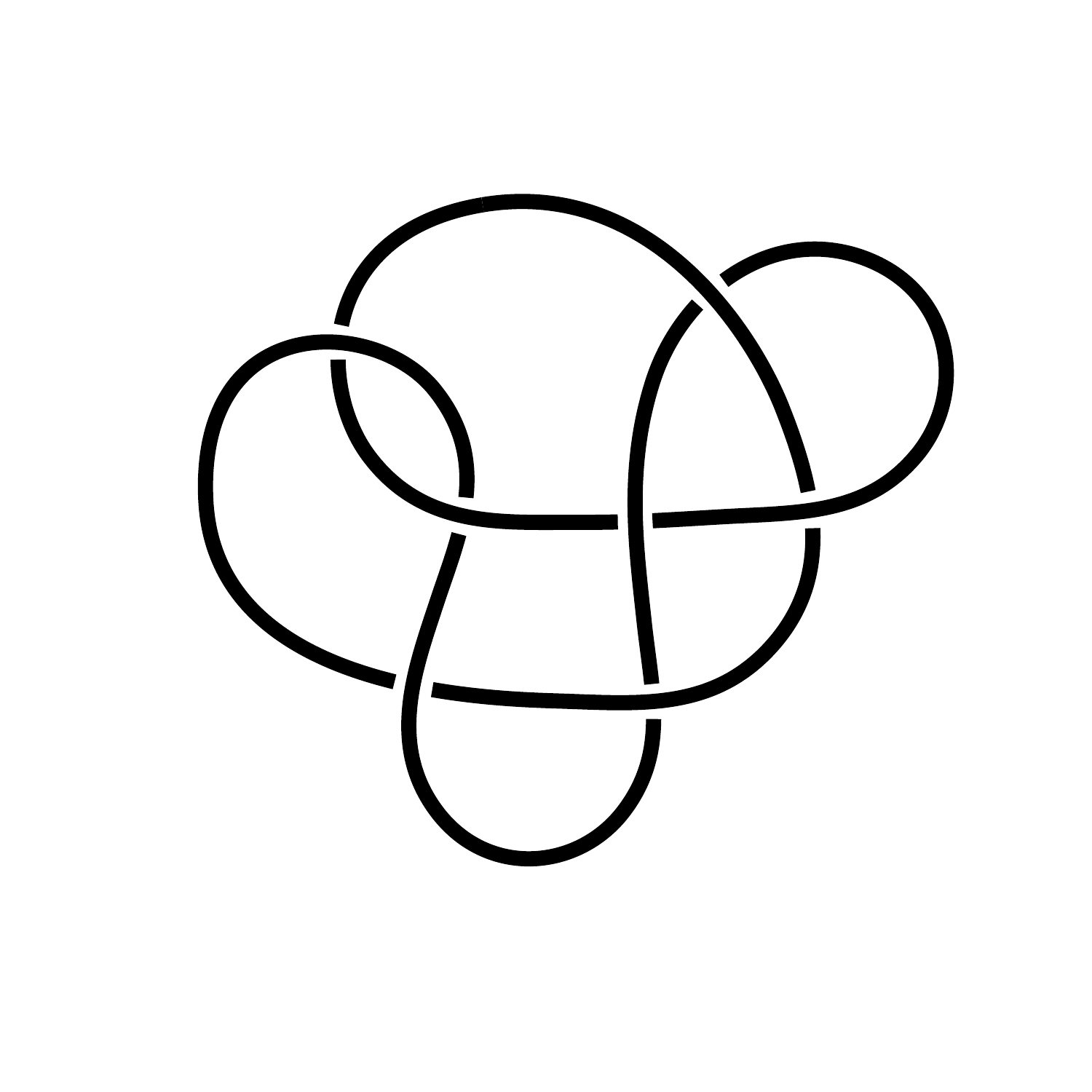}
\caption{Knot $7_6$, post-optimization with breaks indicating
  \label{76_overunder}  underpassing strands}
\end{figure}

With reference to Figure~\ref{k76_strands}, each row of {\tt ou76}
corresponds to a crossing; the first element gives the overstrand and
the second the understrand; thus strand~12 passes over strand~1,
strand~2 passes over strand~11, and so on. The result is shown in
figure~\ref{76_overunder}.

\subsection{Symmetry}

Many of the knots in Figure~\ref{knot_table} have an axis of symmetry,
or possess rotational symmetry.  The package has the ability to impose
symmetry relations on knots, and to optimize the resulting symmetrical
knot.  Minimizing the badness is not entirely straightforward on
account of the induced redundancy, which is characterized using a
symmetry object specific to the knot under consideration.  However,
symmetry constraints do reduce the dimensionality of the optimization
problem.

\begin{figure}[!tbp]
  \centering
\includegraphics[width=11cm]{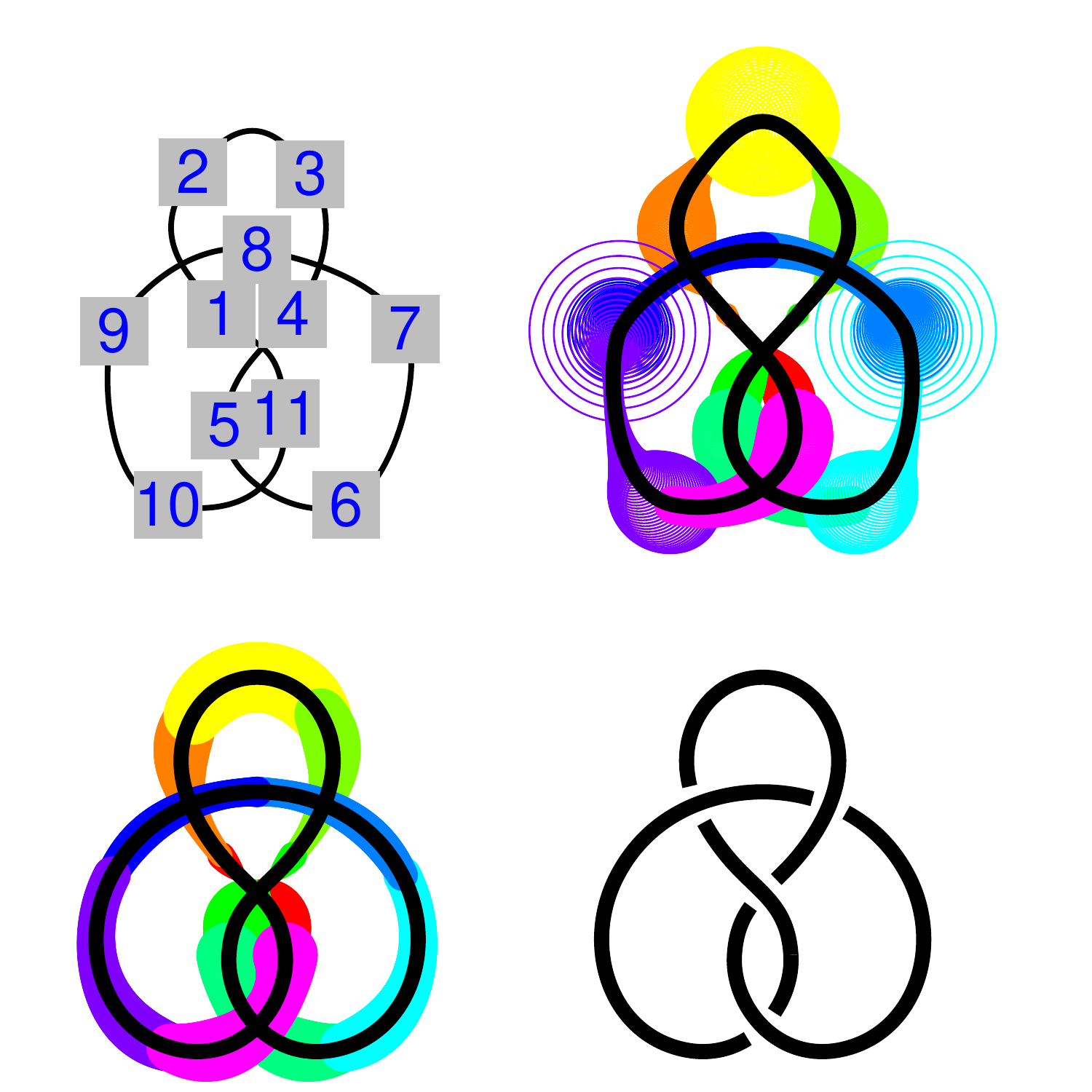}
\caption{Figure eight knots drawn using different plotting methods.
  Top left, knot path with node numbers shown in order to facilitate
    \label{four_figure_8_knots} definition of the symmetry object; top
    right, result of symmetrizing the rough path; lower left, the
    optimized knot with imposed vertical symmetry, with curvature
    plotted; lower right, knot plotted with overstrand and understrand
    indicated using line breaks}
\end{figure}

I will consider the figure-of-eight knot~$4_1$ as an example.  Using
Figure~\ref{four_figure_8_knots}, top left, as reference, the
appropriate symmetry object is defined as follows:

\begin{blockcode}
> fig8 <- reader(system.file("4_1_first_draft.svg",package="knotR"))
> Mver8 <- matrix(c(
+     02,03,
+     09,07,
+     05,11,
+     10,06
+     ),ncol=2,byrow=TRUE)
> sym8 <- symmetry_object(fig8, Mver=Mver8, xver=8)
\end{blockcode}

(Matrix {\tt Mver8} specifies that nodes~2 and~3 are symmetric, as
are nodes~9 and~7, and so on; {\tt xver=8} forces node~8 to be on the
axis of symmetry).  The results are shown in
Figure~\ref{four_figure_8_knots}.

\subsection{Rotational symmetry}

Consider knot~$5_1$. This knot has fivefold rotational symmetry, in
addition to a vertical line of symmetry.  The package includes
functionality to impose appropriate symmetry constraints. Using
Figure~\ref{k5_1_twoknots} as reference, we have:

\begin{figure}[htbp]
  \begin{center}
\begin{figurecode}
> knotplot2(k5_1,node=TRUE,width=FALSE)
\end{figurecode}
\includegraphics[width=11cm]{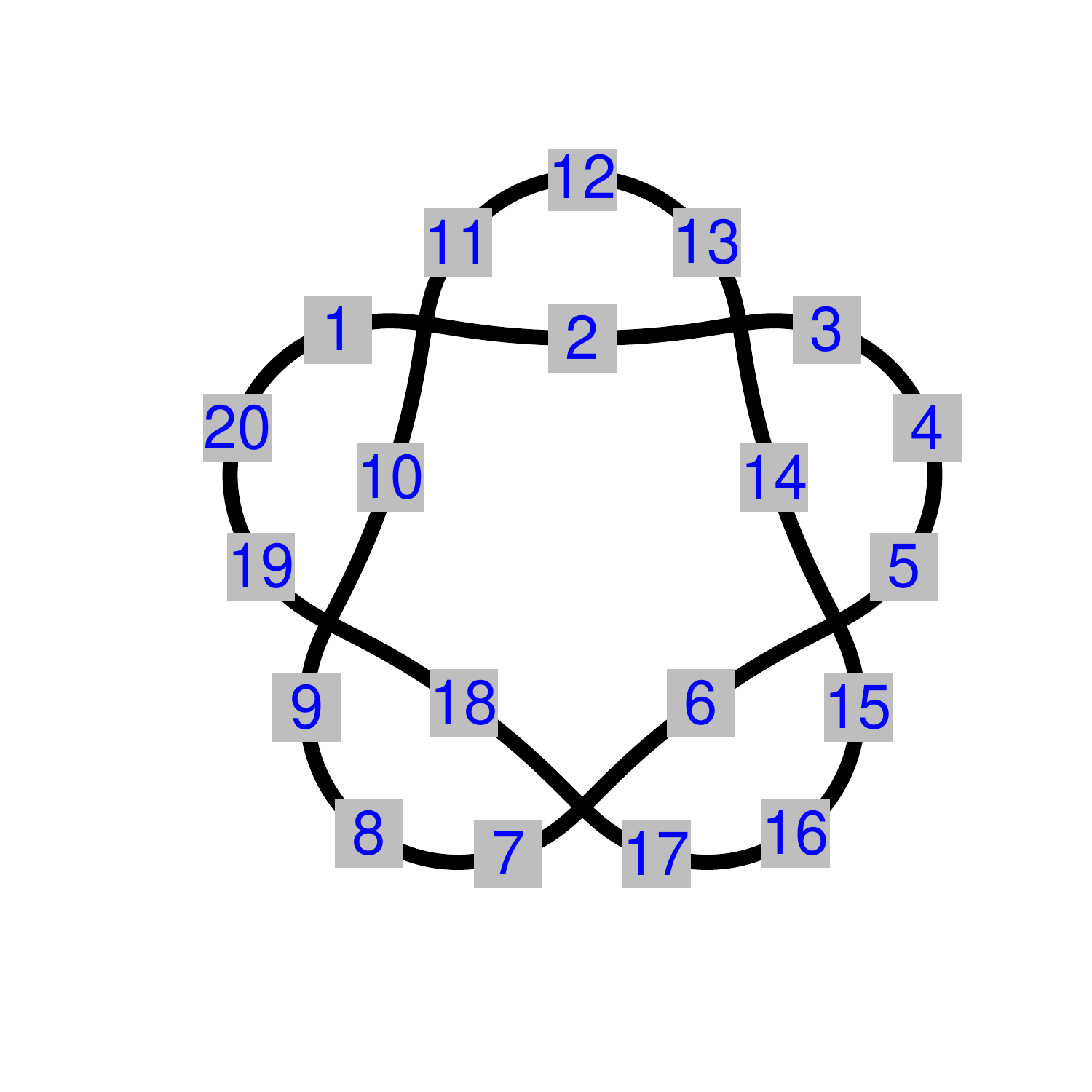}
\caption{Knot $5_1$ \label{k5_1_twoknots} shown with node numbers}
  \end{center}
\end{figure}

\begin{blockcode}
> sym51 <- symmetry_object(k5_1,
+                          Mver = cbind(11,13),
+                          xver = c(2,12),
+                          Mrot = rbind(
+                              c(12,04,16,08,20),
+                              c(13,05,17,09,01),
+                              c(11,03,15,07,19),
+                              c(02,14,06,18,10)
+                          ))
\end{blockcode}

Thus, using the same notation as before, nodes~11 and~13 are
symmetrical about the vertical axis, nodes~2 and~12 are on the
vertical axis.  The {\tt Mrot} argument specifies sets of nodes that
map to themselves under rotation.  The top line of {\tt Mrot}
indicates that nodes 12, 4, 16, 8, and~20 are concyclic.  An example
of a rotationally symmetric knot is given in Figure~\ref{orn20}.

\section{Conclusions and further work}

The {\tt knotR} package allows the user to create rough diagrams of
knots using the inkscape suite of software, and subsequently polish up
such diagrams in terms of a customizable objective function using
numerical optimization techniques.  Further work might include
functionality to deal with links.

\clearpage
\section{Gallery}

There now follows a selection of pleasing knot diagrams taken from
datasets provided with the package.

\begin{figure}[htbp]
  \begin{center}
\begin{figurecode}
par(mfcol=1:2)
knotplot(perko_A)
knotplot(perko_B)
\end{figurecode}
\includegraphics[width=11cm]{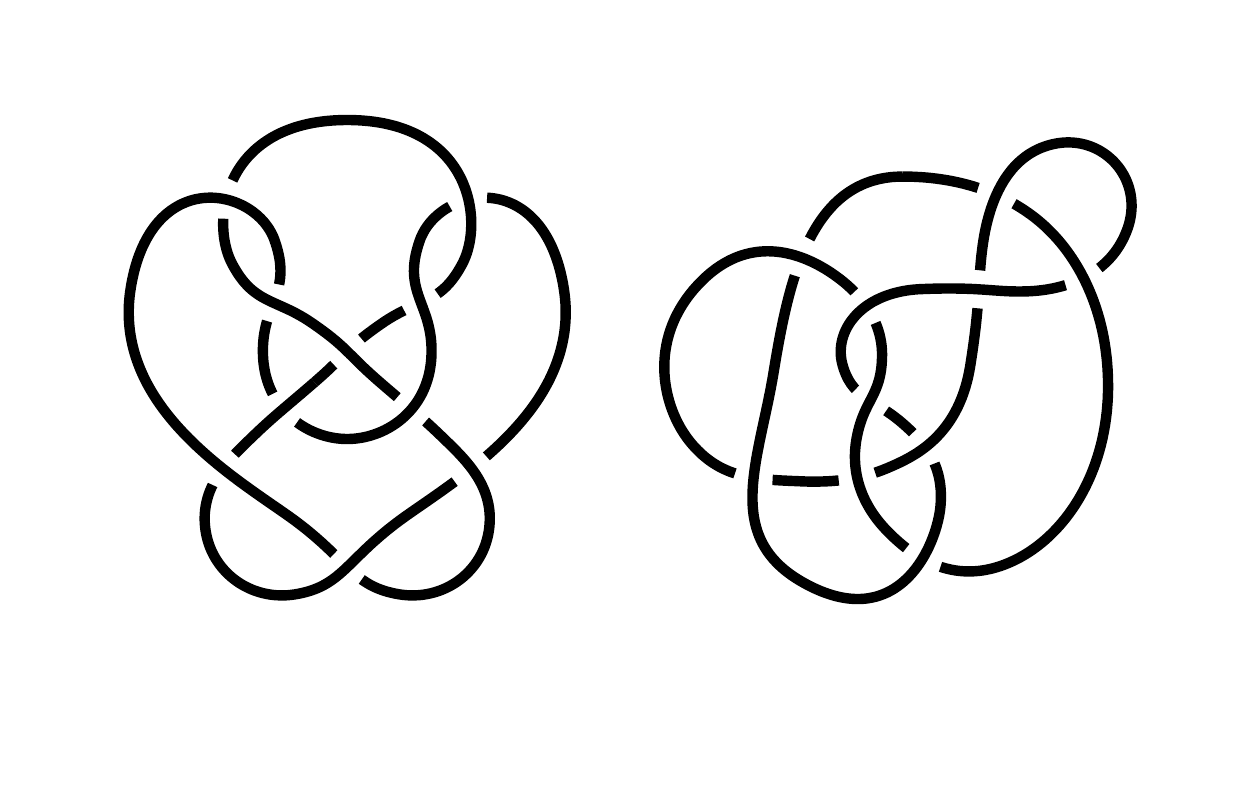}
\caption{Two representations of knot~$10_{125}$, known as the 
\label{perko_AB}  Perko Pair}
  \end{center}
\end{figure}

\begin{figure}[htbp]
  \begin{center}
\begin{figurecode}
> knotplot(ornamental20,gap=15)
\end{figurecode}
\includegraphics[width=11cm]{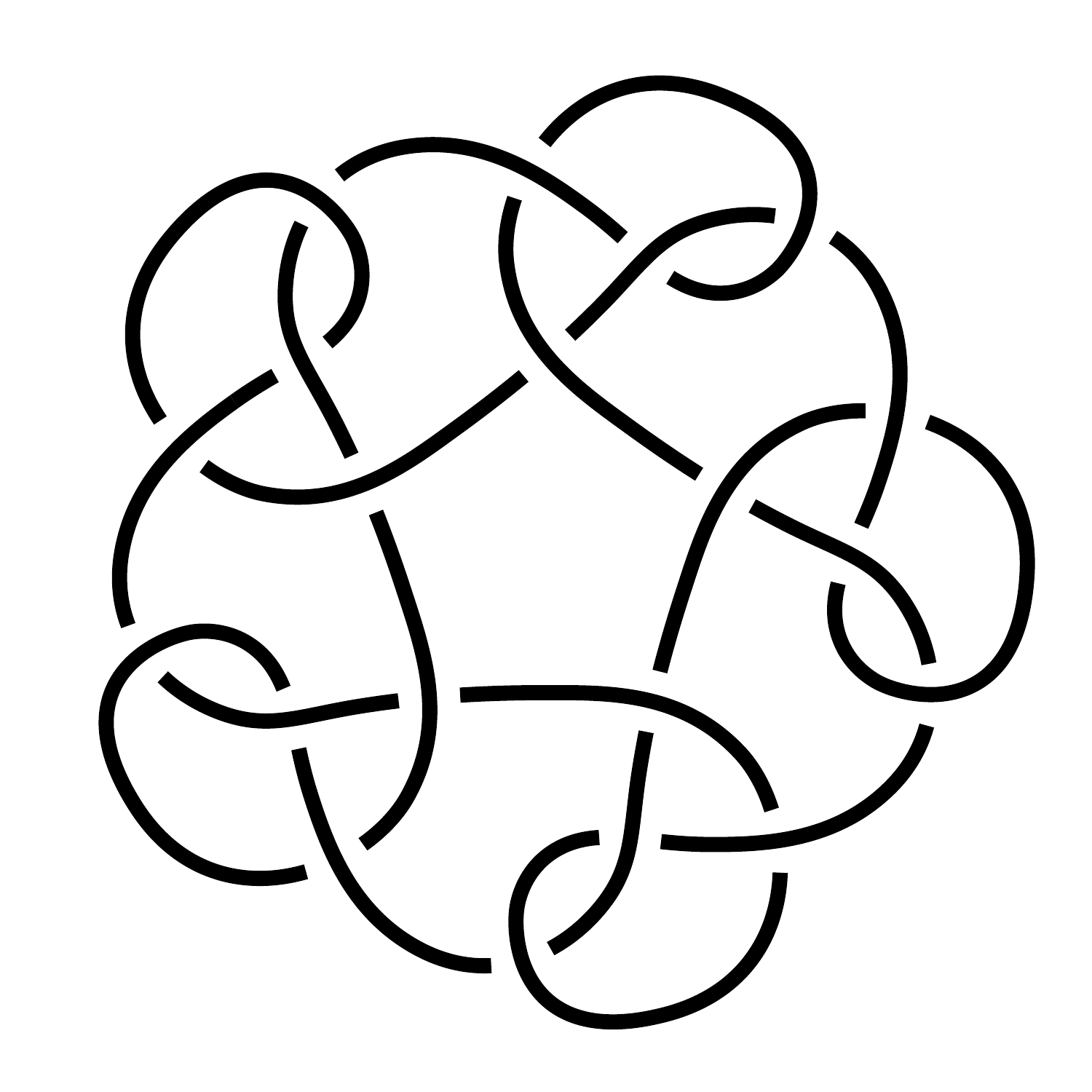}
\caption{An ornamental knot exhibiting fivefold
  rotational \label{orn20} symmetry; note the absence of mirror symmetry}
  \end{center}
\end{figure}

\begin{figure}[htbp]
  \begin{center}
\includegraphics[width=13cm]{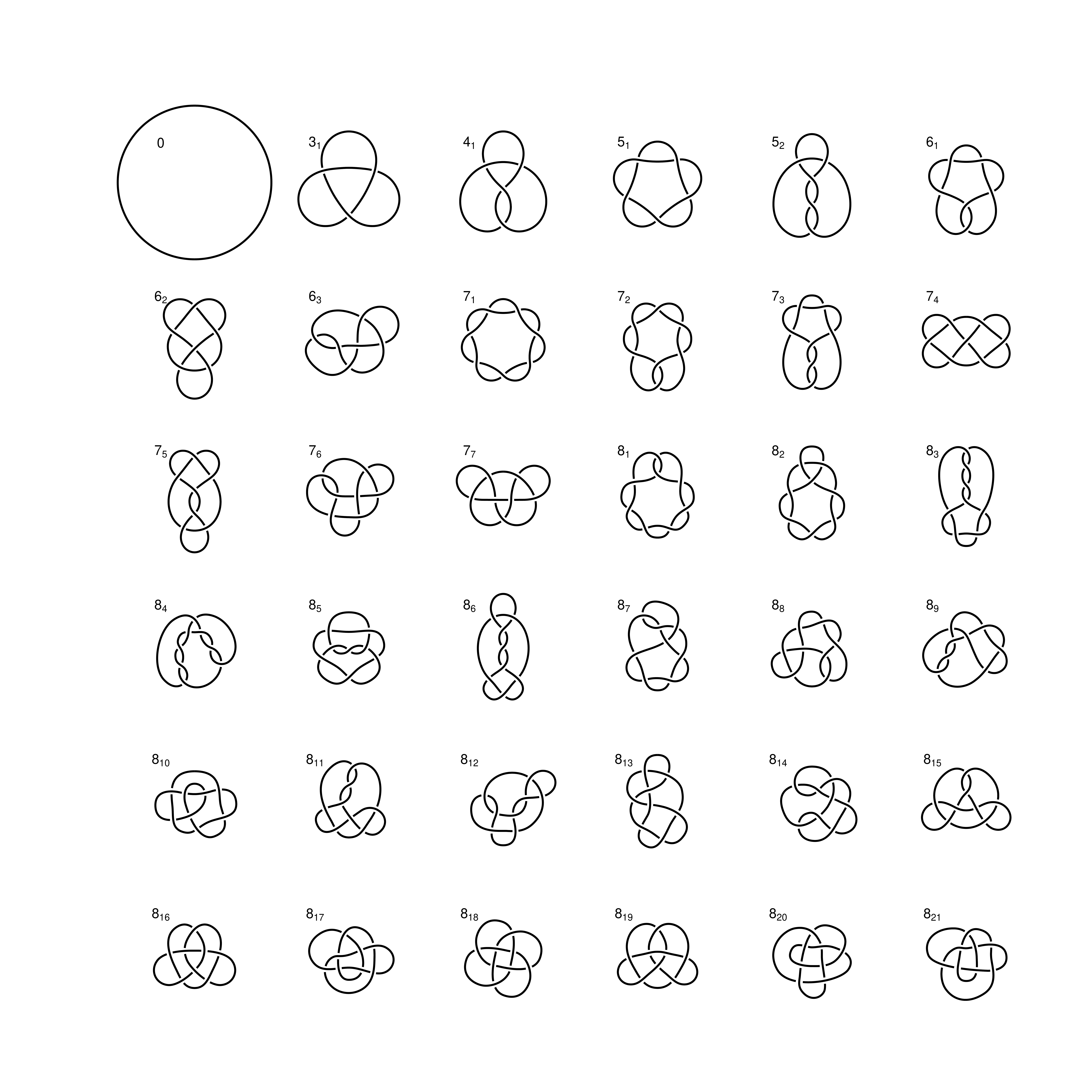}
\caption{All \label{all8} prime knots with 8 or fewer crossings,
notation following Rolfsen}
  \end{center}
\end{figure}

\bibliography{hankin_knot}

\begin{thebibliography}{1}

\bibitem{alexander1926}
J.~W. Alexander and G.~B. Briggs.
\newblock On {Types} of {Knotted} {Curves}.
\newblock {\em Annals of Mathematics}, 28(1/4):562--586, 1926.

\bibitem{hankin2016}
Robin K.~S. Hankin.
\newblock {\em knotR: Knot Diagrams using Bezier Curves}, 2016.
\newblock R package version 1.0-1.

\bibitem{kirsanov2009}
Dmitry Kirsanov.
\newblock {\em The {Book} of {Inkscape}: {The} {Definitive} {Guide} to the
  {Free} {Graphics} {Editor}}.
\newblock No Starch Press, San Francisco, 1 edition edition, October 2009.

\bibitem{manturov2004}
V.~Manturov.
\newblock {\em Knot Theory}.
\newblock Chapman and Hall, 2004.

\bibitem{olsen2014}
Aaron Olsen.
\newblock {\em bezier: Bezier Curve and Spline Toolkit}, 2014.
\newblock R package version 1.1.

\bibitem{rcore2016}
{R Core Team}.
\newblock {\em R: A Language and Environment for Statistical Computing}.
\newblock R Foundation for Statistical Computing, Vienna, Austria, 2016.

\bibitem{wikipedia_knot_theory}
Wikipedia.
\newblock Knot theory --- wikipedia{,} the free encyclopedia, 2016.
\newblock [Online; accessed 17-June-2016].

\end{thebibliography}
\address{Robin K. S. Hankin
Auckland University of Technology\\
2-14 Wakefield Street\\
Auckland\\
New Zealand}
\email{hankin.robin@gmail.com}

\bibliographystyle{plain}

\end{document}